%March 17, 2003
\input amstex
\documentstyle{amsppt}
\magnification=1200 \hsize=13.8cm \catcode`\@=11
\def\NoLogo{\let\logo@\empty}
\catcode`\@=\active \NoLogo

\def\lf{\left}
\def\ri{\right}

\def\a{\alpha}

\def\e{\epsilon}
\def\p{\partial}

\def\C{\Bbb C}
\def\P{\Bbb P}
\def\R{\Bbb R}
\def\oR{\overset\text{o}\to R}

\def\abb{{\alpha\bar\beta}}

\def\tg{\tilde g}
\def\bg{\bar g}
\def \D {\Delta}
\def\aint{\frac{\ \ }{\ \ }{\hskip -0.4cm}\int}
\documentstyle{amsppt}
\magnification=1200 \hsize=13.8cm \vsize=19 cm

\leftheadtext{Lei Ni} \rightheadtext{Ancient solutions} \topmatter
\title{Ancient solutions to
K\"ahler-Ricci flow}\endtitle

\author{Lei Ni\footnotemark }\endauthor
\footnotetext"$^{1}$"{Research partially supported by NSF grants
 and an Alfred P. Sloan Fellowship, USA.}
\address
Department of Mathematics, University of California, San Diego, La
Jolla, CA 92093
\endaddress
\email{ lni\@math.ucsd.edu}
\endemail

\affil { University of California, San Diego}
\endaffil

\date  August 2004 (Revised, March 2005) \enddate

\abstract In this paper, we prove that any non-flat ancient
solution to K\"ahler-Ricci flow with bounded nonnegative
bisectional curvature has asymptotic volume ratio zero. We also
prove that any gradient shrinking solitons with positive
bisectional curvature must be compact. Both results generalize the
corresponding earlier results of Perelman in \cite{P1} and
\cite{P2}. The results can be applied to study the geometry and
function theory of complete K\"ahler manifolds with nonnegative
bisectional curvature via K\"ahler-Ricci flow. It also implies a
compactness theorem on ancient solutions to K\"ahler-Ricci flow.
\endabstract

\endtopmatter

\document

\subheading{\S0 Introduction}\vskip .2cm

The K\"ahler-Ricci flow
$$
\frac{\p}{\p t}g_{\abb}(x,t)=-R_{\abb}(x,t) \tag 0.1
$$
has been  useful in the study of complex geometry in the work of
\cite{B, C1, M2}, etc. The ancient solutions arise \cite{H5} when
one applies the parabolic blow-up to a  finite time singularity or
slowly forming (Type II) singularity as $t$ approaches infinity.
In \cite{H5} Hamilton  introduced some geometric invariants
associated with ancient solutions. One of them is the so-called
asymptotic volume ratio (also called cone angle at infinity),
which is defined  as ${\Cal V}(M, g):=\lim_{r\to
\infty}\frac{V_o(r)}{\omega_n r^n}$, for a complete Riemannian
manifold $(M, g)$. Here $n$ is the dimension (real) of $M$,
$V_o(r)$ is the volume of the ball of radius $r$ centered at $o$
(with respect to metric $g$) and $\omega_n$ is the volume of unit
ball in $\R^n$. This asymptotic volume ratio is well-defined and
independent of the choice of $o$ in the case when $M$ has
non-negative Ricci curvature. One should consult \cite{H5} for the
condition (and the proof) under which  ${\Cal V}(M, g(t))$ is
independent of $t$ for a family of metrics $g(t)$ satisfying Ricci
flow. When the meaning is clear in the context we simply denote
${\Cal V}(M, g(t))$ by ${\Cal V}(g(t))$.

In \cite{P1}, Perelman studied the properties of ancient solutions
with {\it bounded nonnegative curvature operator}, via his entropy
and reduced volume (reduced distance) monotonicity. In particular,
  the following result is proved.

\proclaim{Theorem 1} Let $(M, g(t))$ be a complete non-flat
ancient solution to Ricci flow, with bounded nonnegative curvature
operator. Then ${\Cal V}(g(t))=0$.
\endproclaim

This result holds the key to the rest striking results on  ancient
solutions in \cite{P1}. The assumption on the nonnegativity of the
curvature operator is ensured in dimension three by
Hamilton-Ivey's pinching estimate \cite{H5, Theorem 24.4}, if the
ancient solution is obtained as a blow-up limit of a finite time
singularity. It is also needed to make effective uses of the
reduced distance introduced in \cite{P1, Section 7} in Perelman's
blow down procedure \cite{P1. Section 11} in the study of the
ancient solutions. For K\"ahler-Ricci flow, one would like to
replace the {\it nonnegativity of the  curvature operator} by the
{\it nonnegativity of bisectional curvature} since  the
nonnegativity of the sectional curvature is neither natural nor
necessarily preserved under  K\"ahler-Ricci flow. On the other
hand, the argument of \cite{P1, Section 11} made essential uses of
the {\it nonnegativity of sectional curvature} and the properties
of such Alexandrov spaces. There is no obvious way of adapting the
proof of \cite{P1, Section 11.4} to the K\"ahler setting.
Therefore one needs some  new ingredients in order to generalize
Perelman's result to K\"ahler-Ricci flow assuming only the
nonnegativity of the bisectional curvature. It turns out that this
technical hurdle can be overcome by a result (Proposition 1.1) on
the shrinking solitons, Perelman's blow-down procedure and a
splitting result for K\"ahler manifolds with nonnegative
bisectional curvature, proved recently by Luen-Fai Tam and the
author in \cite{NT2}. These will be the main focus of this paper.

There exists another motivation of proving Theorem 1 for
K\"ahler-Ricci flow on K\"ahler manifolds with  nonnegative
bisectional curvature. In a  AMS meeting (November 2001, held at
Irvine, California), Huai-dong Cao proposed the following
conjecture (see also \cite{C4} for a related problem on ancient
solutions).

\proclaim{Conjecture (Cao)} Let $(M^m, g(t))$ be a non-flat steady
gradient  K\"ahler-Ricci soliton with bounded nonnegative
bisectional curvature. Then ${\Cal V}(g(t))=0$.
\endproclaim

Proving Theorem 1 for K\"ahler-Ricci flow on  K\"ahler manifolds
with nonnegative bisectional curvature will imply  Cao's
conjecture. Combining the techniques from \cite{P1} and \cite{P2}
with \cite{NT2}, we indeed can prove the  following such
generalization  of Theorem 1.

\proclaim{Theorem 2} Let $(M^m, g(t))$ ($m=\dim_{\C} (M)$, n=2m)
be a non-flat ancient solution to K\"ahler-Ricci flow {\rm (0.1)}.
Assume that $(M, g(t))$ has bounded nonnegative bisectional
curvature. Then  ${\Cal V}(g(t))=0$.
\endproclaim

 The only known result along this line, under
 the assumption of {\it nonnegative bisectional curvature} has been only proved  for
 the case of $m=2$,  and only for the
gradient K\"ahler-Ricci solitons, which are special ancient
solutions. Please see, for example \cite{C4, CZ2}. Following
\cite{P1}, an immediate consequence of Theorem 2 is the following
compactness result.

\medskip

{\it For an fixed $\kappa>0$, the set of $\kappa$-solution to
K\"ahler-Ricci flow is compact module scaling. }

\medskip

Please refer to Section 2 for the definitions of the
$\kappa$-solutions. A compactness result for the elliptic case,
namely the K\"ahler-Einstein metrics, was proved earlier in
\cite{T} for  compact manifolds under extra assumptions on an
integral bound of curvature, a volume lower bound and a diameter
upper bound.

 Since for the Riemannian manifolds with
non-negative Ricci curvature, the Bishop volume comparison theorem
asserts that $\frac{V_o(r)}{\omega_{2m}r^{2m}}$ is monotone
non-increasing,   the manifold  $M$ is called of maximum volume
growth if ${\Cal V}(M, g)>0$. The above result simply concludes
that the ancient solutions with bounded nonnegative bisectional
curvature (which is preserved under the K\"ahler-Ricci flow, by
\cite{B, M2}) can not be of maximum volume growth. Applying Shi's
short time existence result \cite{Sh1}, Hamilton's singularity
analysis argument in \cite{H5, Theorem 16.2} (or Perelman's result
in \cite{P1, Section 11}), as well as estimates from  \cite{NT1}
and \cite{N2} one can have the following corollary of Theorem 2.

\proclaim{Corollary 1} Let $(M^m, g_0)$ ($m=\dim_\C(M)$) be a
complete  K\"ahler manifold with bounded nonnegative bisectional
curvature. Assume that $M$ is of maximum volume growth. Then the
K\"ahler-Ricci flow {\rm (0.1)}, with $g(x,0)=g_0(x)$ has a long
time solution. Moreover the solution has no slowly forming (Type
II) singularity as $t$ approaches $\infty$. In particular,
 there exists a $C=C(M, g_0)>0$ such that the scalar
curvature ${\Cal R}(y)$ satisfies
$$
\aint_{B_x(r)}{\Cal R}(y)\, d\mu \le \frac{C}{(1+r)^2}, \tag 0.2
$$
where $B_x(r)$ is the ball of radius $r$ centered at $x$,
$V_x(r)=\text{Vol} (B_x(r))$ and $\aint_{B_x(r)}f(y)\,
d\mu=\frac{1}{V_x(r)}\int_{B_x(r)}f(y)\, d\mu.$ Furthermore, $M$
is differmorphic (homeomorphic) to $\C^m$, for $m> 2$ ($m=2$), and
is biholomorphic to a pseudoconvex domain in $\C^m$.
\endproclaim

The curvature decay  statement of Corollary 1 confirms a
conjecture of Yau in \cite{Y, page 621}, where he speculated that
two assumptions in Shi's main theorem of \cite{Sh3} can be
replaced by the maximum volume growth alone. See also the recent
work of Wu and Zheng \cite{WuZ} on various examples related to the
above corollary. Corollary 2 also provides a partial answer to a
question asked in \cite{N2, Conjecture 3.1} on conditions
equivalent to the existence of holomorphic functions of polynomial
growth. In fact, there the author speculated that either the
maximum volume growth or the average quadratic curvature decay is
equivalent to the existence of nonconstant polynomial growth
holomorphic functions provided the manifold has quasi-positive
bisectional curvature. There also exists a related general
conjecture of Yau on the non-existence of the bounded holomorphic
functions, on which one can refer to \cite{LW} for some recent
progresses. Using Perelman's Theorem 1 one can conclude the
similar result for the Riemannian manifolds with {\it bounded
nonnegative curvature operator} and maximum volume growth. For
K\"ahler-Ricci flow, under the stronger assumption that $M$ is a
K\"ahler manifold with {\it bounded nonnegative curvature
operator} and maximum volume growth, a similar (but slightly
weaker) estimate as (0.2) was proved earlier in \cite{CZ4} by
applying
 the dimension reduction argument of Hamilton. This  is in fact a
 fairly easy consequence of
 Perelman's Theorem 1. (See
the proof of Corollary 1 for details.  The main constraint of
related results in \cite{H5}, for the application to K\"ahler
geometry, is that the argument of \cite{H5} only works under the
stronger assumption on {\it nonnegativity of curvature
operator/sectional curvature}, which is sufficient for the   study
of three manifolds, but not for the  study of K\"ahler manifolds
of complex dimension $\ge 3$.) In \cite{Sh3-4}, the long time
existence result was proved under a uniform curvature decay
assumption similar to (0.2), with/without the assumption of the
maximum volume growth. In Corollary 1  we have the curvature decay
as a consequence instead of an assumption, exactly as what Yau
speculated in \cite{Y}. Under certain  further curvature average
decay assumptions,  the topological conclusion in Corollary 1 was
proved earlier in \cite{CZ3}. (See also \cite{H6, Sh4} for related
earlier fundamental works. The key observation on the improvement
of injectivity radius lower bound was first observed in
\cite{Sh4}.) The last statement of Corollary 1 was proved first in
\cite{Sh4} (see also \cite{CZ3}) under the assumptions that the
manifold has positive sectional curvature and a certain curvature
average decay condition. Note also that Corollary 4.1 of
\cite{NT2} proved that $M$ is Stein if $M$ has maximum volume
growth and nonnegative bisectional curvature (which is a
conjecture of H. Wu).

\proclaim{Corollary 2} Let $(M^m, g_0)$ be a complete K\"ahler
manifold with bounded nonnegative holomorphic bisectional
curvature and maximum volume growth. Then the transcendence degree
of the rational function field of $M$ (see \cite{N2} for
definition) is equal to $m$. Furthermore, $M$ is biholomorphic to
an affine quasi-algebraic variety. When $m=2$, $M$ is
biholomorphic to $\C^2$.
\endproclaim

  In \cite{M1, M3}, a systematic scheme on embedding/compactifying
  complete K\"ahler manifolds with
  positive curvature was developed.  One can also consult the survey article \cite{M4}
  for expositions on these methods and  related results.
  Originally in \cite{M1}  under the  assumption
that $M$ has {\it positive bisectional curvature} and certain
curvature  decay conditions, the conclusion that $M$ is
biholomorphic to $\C^2$ in the case $m=2$, and that $M$ is
biholomorphic to an affine variety in the higher dimension  were
obtained. In fact in his fundamental work \cite{M1} Mok laid down
all the framework of the affine embedding and observed that one
can apply the result of Ramanujam, which asserts that a
quasi-projective surface homeomorphic to $\R^4$ must be
biholomorphic to $\C^2$, in the case $m=2$. Later, for the case
$m=2$ only, in \cite{CTZ}, following Mok's compactification scheme
in \cite{M1. M3}, the authors improved the above result of Mok,
replacing the assumptions of \cite{M1} by that $M$ has {\it
bounded positive bisectional curvature} and of maximum volume
growth. The result stated in \cite{CTZ} also assumes a mild
average curvature decay condition which   replaces the stronger
point-wise curvature decay assumption of \cite{M1} via
K\"ahler-Ricci flow. This average curvature decay condition can be
removed  by combining with another later paper \cite{CZ4}, again
assuming the {\it positivity} of the bisectional curvature and
volume being of maximum growth. In this later improvement
\cite{CTZ}, which is only restricted to the case of $m=2$ (also
the earlier paper \cite{CZ2}), the proof also crucially relies on
an observation only valid in complex dimension two, originally due
to Ivey \cite{I}, that an ancient solution to K\"ahler-Ricci flow
with nonnegative bisectional curvature must have nonnegative
curvature operator. Namely the method there relied crucially on
the dimensional reduction results of Hamilton in \cite{H5}. Hence
it does not work under the assumption of the nonnegativity of
bisectional curvature in higher dimensions.

In Corollary 2, when $m=2$, our statement  assumes only  {\it
nonnegativity} instead of {\it positivity} of the bisectional
curvature. Moreover, our new approach also works for the  higher
dimensional case. Note that the uniform multiplicity estimates
recently proved in \cite{N2}  simplifies the steps in \cite{M1}
quite a bit.  We should point out that the $m=2$ case can also be
obtained by combining Corollary 4.1 in \cite{NT2} obtained by
Luen-Fai Tam and the author  with the argument of \cite{CTZ}.
However the method of current paper provides an unified direct
approach which works for any dimensions.

In the proof of Theorem 2 we need the following  result on
gradient shrinking solitons of K\"ahler-Ricci flow, which is of
independent interests.

\proclaim{Theorem 3} Let $(M^m, g)$ be a non-flat  gradient
shrinking soliton to K\"ahler-Ricci flow. \roster \item"{(i)}" If
the bisectional curvature of $M$ is  positive then $M$ must be
compact and isometric-biholomorphic to $\P^m$. \item"{(ii)}" If
$M$ has nonnegative bisectional curvature then the universal cover
$\tilde M$ splits as $\tilde M= N_1\times N_2\times \cdots \times
N_l\times \C^k$ isometric-biholomorphically, where $N_i$ are
compact irreducible Hermitian Symmetric Spaces.
\endroster
In particular, ${\Cal V}(M,g)=0$.
\endproclaim

Theorem 3 generalizes a corresponding recent result of \cite{P2,
Lemma 1.2}, where Perelman shows that in dimension three, any
$\kappa$-noncollapsed  gradient shrinking soliton with {\it
bounded positive sectional curvature} must be compact.  When $M$
is compact and $m=2$, the classification in part (ii) was obtained
in \cite{I} under even weaker assumption on the nonnegativity of
the isotropic curvature.

For recent works on K\"ahler-Ricci flow on {\it compact}
manifolds, one should refer to the survey articles \cite{CC, Cn}.
 Perelman \cite{P3} also has some
important  work on the conjecture (which arises from the related
works of Hamilton and Tian) concerning the large time behavior of
(normalized) K\"ahler-Ricci flow on compact K\"ahler manifolds
with $c_1(M)>0$ (cf. \cite{N3}). In \cite{CT1-2}, very recent
progresses have been made towards the uniformization problem
addressed in Corollary 2.

Finally, we should point out that when  dimension $m=1$ the above
results are  known from the work of Hamilton \cite{H3} and
\cite{H5, Section 26}. (See also the work of  Chow \cite{Ch} and
books \cite{CK, CLN}.)

The method of this paper has other applications in the study of
Ricci flow on Riemannian manifolds. Please see Section 3.

\medskip

{\it Acknowledgement.} The author would like to thank Professors
H.-D. Cao, B. Chow, T. Ilmanen, R. Schoen and Jon Wolfson for
helpful discussions, Professor P. Lu for pointing out the
reference \cite{Ye}, Professors L.-F. Tam,  M.-T. Wang, H. Wu and
F.-Y. Zheng
 for their interests. Professor H.-D. Cao pointed out to the author
that the earlier version of Corollary 2 has further complex
geometric consequences relating to the uniformization of K\"ahler
manifolds with nonnegative bisectional curvature.  It is a
pleasure to record our gratitude to him. Finally we would like to
thank the referee for pointing out a discrepancy in the earlier
version of this paper (of which the author was also aware
independently).

\subheading{\S1 Proof to Theorem 2 and 3}\vskip .2cm

Recall that a complete Riemannian manifold $(M, g)$ is called a
gradient shrinking soliton if there exists a smooth function $f$
such that, for some positive constant $a$,
$$
\nabla_i\nabla_j f+R_{ij}-ag_{ij}=0. \tag 1.1
$$

\proclaim{Proposition 1.1} Let $(M,g)$ be a Ricci non-flat
gradient shrinking soliton. Assume that the Ricci curvature of $M$
is nonnegative.  Then there exists a $\delta=\delta(M)\, (1\ge
\delta
>0$) such that
$$
{\Cal R}(x)\ge \delta>0. \tag 1.2
$$
\endproclaim
\demo{Proof} It is well know from the strong maximum principle
that the scalar curvature ${\Cal R}(x)>0$. Differentiating (1.1)
and applying the second Bianchi identity, we have that
$$
\nabla_i {\Cal R}=2R_{ij}f_j. \tag 1.3
$$
This implies that
$$
\split \nabla_i{\Cal
R}+2f_{ij}f_j-2af_i&=2\left(R_{ij}+f_{ij}-ag_{ij}\right)f_j\\
&=0
\endsplit
$$
which further implies that there exists a constant $C_1=C-1(M)$
such that
$${\Cal R}+|\nabla
f|^2-2af=C_1. \tag 1.4
$$
These equations are well known for the gradient shrinking
solitons.

 Let $o\in M$ be a fixed point. For any $x\in M$ we denote the distance of $x$
 from $o$ by $r(x)$. Let $\gamma(s)$ be minimal geodesic joining
 $x$ from $o$, parametrized by the arc-length. For simplicity we often also denote $r(x)$
 by $s_0$.  Let $\{E_i(s)\}$ ($0\le i\le n-1$)  be a parallel frame along $\gamma(s)$ such that
$E_0(s)=\gamma'(s)$. If $s_0\ge 2$, for  $s_0\ge r_0\ge 1$,
choose $n-1$-variational vector fields $Y_i(s)$ ($1\le i\le n-1$)
along $\gamma(s)$ as
$$
Y_i(s)=\cases & sE_i(s), \quad  0\le s\le 1\\
              & E_i(s), \quad 1\le s\le s_0-r_0\\
              & \frac{s_0-s}{r_0}E_i(s), s_0-r_0\le s\le s_0
              \endcases.
$$
By the second variation consideration
 \cite{P1, Lemma 8.3 (b)} (see also \cite{H5, Theorem 17.4}), we
  have that
$$
 \sum_{i=1}^{n-1} \int_0^{s_0} |Y'_i(s)|^2-R(\gamma'(s), Y_i(s),
\gamma'(s), Y_i(s))\, ds\ge 0.
$$
In particular we can find $C(M)$, which depends only the upper
bound of the Ricci curvature of $M$ on $B_o(1)$,  such that
$$
\split
 \int_0^{s_0-r_0}Ric(\gamma'(s), \gamma'(s))\, ds &\le
C(M)+\frac{n-1}{r_0}-\int_{s_0-r_0}^{s_0}\left(\frac{s_0-s}{r_0}\right)^2 Ric(\gamma'(s), \gamma'(s))\, ds\\
&\le C(M)+\frac{n-1}{r_0}.
\endsplit\tag 1.5
$$
Here we have used the fact the Ricci curvature is nonnegative. We
claim that  there exists a positive  constant $A=A(M)$, if $s_0\ge
A$ and ${\Cal R}(x)\le 1$, there exists another constant, still
denoted by  $C(M)$,  such that
$$
\int_0^{s_0} Ric(\gamma'(s), \gamma'(s))\, ds \le
\frac{a}{2}s_0+C(M). \tag 1.6
$$
Assume that we have proved the claim (1.6). Then there exists
$C_2=C_2(M)>0$ such that
$$
\split
 \langle \nabla f,
 \gamma'(s)\rangle|^{\gamma(r(x))}_o&=\int_0^{s_0}\frac{d}{ds}(\langle \nabla f,
 \gamma'(s)\rangle )\, ds\\
 &=\int_0^{s_0}(\nabla_i\nabla_j f)\frac{d \gamma^i(s)}{d s}\frac{d \gamma^j(s)}{d
 s}\, ds\\
 &=\int_0^{s_0}(a-Ric(\gamma'(s), \gamma'(s)))\, ds\\
 &=ar(x)-\int_0^{s_0}Ric(\gamma'(s), \gamma'(s))\, ds\\
 &\ge \frac{a}{2} r(x)-C_2,
 \endsplit
$$
which implies that for  every $x \in M\setminus B_o(A)$ with
${\Cal R}(x)\le 1$,
$$
\langle \nabla f, \nabla r\rangle (x) \ge
\frac{a}{2}r(x)-C_2-|\nabla f|(o).
$$
It in particular implies that for every such $x$, with $r(x)\ge
\frac{4}{a}(C_2+|\nabla f|(o))$ additionally, $\nabla f\ne 0$.

Now we can prove the proposition after the claim (1.6). For any
$x\in M\setminus (B_o(A)\cup B_o(\frac{4}{a}(C_2+|\nabla
f|(o))))$, without the loss of generality we can assume that
${\Cal R}(x)\le 1$,  let $\sigma(\eta)$ be the integral curves of
$\nabla f$, passing $x$ with $\sigma(0)=x$. By (1.3) we have that
$$
-\frac{d}{d\eta}({\Cal R}(\sigma(\eta))=-2R_{ij}\frac{d
\sigma^{i}}{d\eta}\frac{d \sigma^{j}}{d\eta}\le 0. \tag 1.7
$$
This implies that ${\Cal R}(x)\ge {\Cal R}(\sigma(\eta))$, for
$\eta <0$. On the other hand,
$$
-\frac{d }{d\eta}r(\sigma(\eta))=-\langle \nabla r, \nabla
f\rangle \le -(C_2+|\nabla f|(o))\le 0 \tag 1.8
$$
as far as $r(\sigma(\eta))\ge\max(A, \frac{4}{a}(C_2+|\nabla
f|(o)))$, noticing that we always have  ${\Cal R}(\sigma(\eta))\le
1$. This implies that the integral curve $\sigma$ exists for all
$\eta<0$ since $|\nabla f|$ is bounded inside the closed ball
$\overline {B_o(2r(x))}$. The estimate (1.8) also implies that
there exists $\eta_1<0$ such that $r(\sigma(\eta_1))=\max(
A,\frac{4}{a}(C_2+|\nabla f|(o)))$. Applying (1.7) we have that
$$
{\Cal R}(x)\ge \inf_{y\in B_o(r(\sigma(\eta_1)))}{\Cal R}(y).
$$
This proves the proposition assuming the claim (1.6).

Now we prove the claim (1.6). First by equation (1.3) and the fact
$R_{ij}\ge 0$ we have that $f_{ij}\le a g_{ij}$. This implies that
along any minimizing geodesic $\gamma(s)$ from $o$, $f''(s)\le a$.
Hence there exists $B=B(M)$ such that
$$
f(x)\le (a+1)r^2(x)
$$
for $r(x)\ge B$. Using (1.4) and the fact that ${\Cal R}> 0$ we
have that
$$
|\nabla f|(x)\le 2(a+1) r(x)
$$
for $r(x)\ge B$. On the other hand, (1.3) also implies that
$$
|\nabla {\Cal R}|^2\le 4 {\Cal R}^2|\nabla f|^2.
$$
The above two inequality implies the the estimate
$$
|\nabla \log {\Cal R}|(x)\le 2(a+1)r(x)
$$
for $r(x)\ge B$. Now we adapt the notations and situations right
before (1.6) and choose $r_0$ in (1.5) such that
$\frac{n-1}{r_0}={\e} s_0$ with  some fixed positive constant
$\e\le \min (1, \frac{a}{2})$. Then (1.6) implies  that
$$
 \int_0^{s_0-r_0}Ric(\gamma'(s), \gamma'(s))\, ds \le C(M)+\e
 s_0.\tag 1.9
$$
Notice that $r_0=\frac{n-1}{\e s_0} \le \frac{n-1}{\e} \le
\frac{s_0}{2}$ if $s_0\ge A$ for some $A=A(M)\ge \max(1, 2B,
2\frac{n-1}{\e})$. Now using the gradient estimate on $\log {\Cal
R}$ above we have that
$$
\split \log\frac{ {\Cal R}(\gamma(s_1))}{ {\Cal
R}(\gamma(s_0))}&=-\int_{s_1}^{s_0}\frac{d}{ds} \log {\Cal
R}(\gamma(s))\, ds\\
& \le \int_{s_1}^{s_0}|\nabla \log {\Cal R}|\, ds\\
&\le 2(a+1)s_0(s_0-s_1)
\endsplit
$$
for $s_1\le s_0$. Hence
$$
{\Cal R}(\gamma(s))\le {\Cal
R}(\gamma(s_0))\exp(\frac{2(a+1)(n-1)}{\e})\le
\exp(\frac{2(a+1)(n-1)}{\e})
$$
for any $s\ge s_0-r_0$. Here we have used the assumption ${\Cal
R}(x)={\Cal R}(\gamma(s_0))\le 1$. This further implies that
$$
\split \int_{s_0-r_0}^{s_0} Ric(\gamma'(s), \gamma'(s))\, ds &\le
\int_{s_0-r_0}^{s_0}{\Cal R}(\gamma(s))\, ds\\
&\le r_0 \exp(\frac{2(a+1)(n-1)}{\e})\\
&=\frac{n-1}{{\e} s_0}\exp(\frac{2(a+1)(n-1)}{\e})\\
&\le C(\e, M).
\endsplit
$$
Together with (1.9), we prove our claim (1.6). Hence we complete
the proof of the Proposition 1.1.
\enddemo

\demo{Proof of Theorem 3}  {\it Case (i).}  By the strong maximum
principle we can assume that ${\Cal R}>0$ everywhere, otherwise
one would have that $M$ is Ricci flat, hence flat. By Proposition
1.1 , we know that ${\Cal R}\ge \delta>0$, for some $\delta$. In
particular, it implies that for any $x$
$$
\aint_{B_x(r)}{\Cal R}(y)\, d\mu \ge \delta.
$$
Here $\aint$ is the average integral defined in Corollary 1 of the
introduction.
 On the other hand, part (ii) of Theorem 4.2 in \cite{NT2}
implies that if $M$ is not compact, then it must satisfies
$$
\aint_{B_x(r)}{\Cal R}(y)\, d\mu\le \frac{C_2}{1+r} \tag 1.10
$$
for some $C_2(x, M)>0$. This is a contradiction! This implies that
$M$ must be compact.

{\it Case (ii)}. Let $\tilde M$ be the universal cover of $M$.
Since $\tilde M$ is still a shrinking soliton,  we can  apply
Proposition 1.1  to $\tilde M$.  From part (i) of Theorem 4.2 in
\cite{NT2}, we can split $\tilde M$ as a product of two manifolds
with one of the factor being compact and the other satisfies the
curvature average decay (1.10). The classification result follows
from the celebrated result of Siu-Yau \cite{SiY}, Mok \cite{M2} as
well as an observation of Koiso \cite{Ko}.

The conclusion ${\Cal V}(M, g)=0$ now follows easily from (i) and
(ii).
\enddemo

The proof of Theorem 2 requires the blow-down procedure of
Perelman in \cite{P1, Proposition 11.2}. Since Proposition 11.2 of
\cite{P1} has only  a sketched proof, in the following we present
in a more detailed way the blow-down procedure of Perelman  of
\cite{P1, 11.2} on the study of so-called {\it bounded
$\kappa$-solutions},  complete ancient solutions with {bounded
nonnegative curvature operator and $\kappa$-non-collapsed
 in all scales}. (See also \cite{CLN}, \cite{KL}, \cite{STW} and \cite{Ye}
 for various expositions on this part.)
 A solution (to Ricci flow) is called $\kappa$-non-collapsed if for any time $t$, for
 any ball $B_x(r)$ (with respect to metric $g(t)$),  satisfying
 $|Rm|(y)\le \frac{1}{r^2}$ for all $y\in
 B_x(r)$ one has $V_x(r)\ge \kappa r^n$.
 One should refer to \cite{P1, Section 4} for more discussions on
  $\kappa$-non-collapsing, \cite{P1, Section 7, 11} for more details of the
 reduced distance and its properties. Now we adapt \cite{P1,
 Section 7, 11} to
 K\"ahler-Ricci flow and replace the nonnegativity of the
 curvature operator by the nonnegativity of the bisectional
 curvature. So in the following,  we let $(M, g(t))$ be a complete
 non-flat ancient solution to K\"ahler-Ricci flow. We also assume
 that $(M, g(t))$ has bounded nonnegative bisectional curvature
 (the boundedness of curvature can be replaced by the differential
 Harnack inequality, (1.14) below)
  and it is $\kappa$-non-collapsed for some $\kappa>0$.

 First we recall the definition of the reduced
 distance $\ell(y, \tau)$. Fix a space-time point $(x_0, t_0)$. Let $\tau=t_0-t$.
 Define
 $$
\ell(y, \tau)=\inf_{\gamma, \gamma(0)=x_0,
\gamma(\tau)=y}\frac{1}{2\sqrt{\tau}}\int_0^\tau
\sqrt{\eta}\lf({\Cal R}+4|\gamma'(\eta)|^2\ri)\, d\eta.
$$
We have factor $4$ here since we study  K\"ahler-Ricci flow
instead of Ricci flow. Since the ${\Cal R}\ge 0$ and the metrics
are shrinking (since $R_{\abb}\ge 0$) along the $t$ direction it
is easy to see that
$$
\ell(y, \tau)\ge \frac{d^2_{t_0}(y, x_0)}{\tau}.
$$
(Sometimes we also denote $d_{t_0}(y, x_0)$ by $d_{0}(y, x_0)$
when we think in terms of $\tau$.) As a consequence of above lower
bound on $\ell$ we know that $\ell(y, \tau)$ achieves its minimum
for each $\tau$ at some point finite distance away from $x_0$.
Thus we can conclude that \proclaim{Claim 1} For each $\tau$ there
exists a point $x(\tau)$ such that $\ell(x(\tau),
\tau)=\frac{n}{2}$.
\endproclaim
\demo{Proof} By the equation (7.15) in \cite{P1}, using the
maximum principle it was shown in \cite{P1, Section 7} that
$\min_{y\in M} \ell(y, \tau)\le \frac{n}{2}$. Using the continuity
and (1.10) we know the existence of $x(\tau)$. Tom Ilmanen has a
nice geometric explanation on this fact by making analogy with the
mean curvature flow.
\enddemo

Let us recall the equations satisfied  by $l(y, \tau)$ from
\cite{P1, Section 7}. First we have
$$
|\nabla \ell|^2+{\Cal
R}=\frac{1}{\tau}\ell-\frac{1}{\tau^{\frac{3}{2}}}K. \tag 1.11
$$
Here
$$
K=\int_0^\tau \eta^{\frac{3}{2}}H(X)\, d\eta
$$
where $H(X)=-{\Cal R}_{\tau}-2\langle\nabla {\Cal R},
X\rangle-2\langle X, \nabla {\Cal R}\rangle +4Ric(X,
X)-\frac{1}{\tau}{\Cal R}$ is the trace different Harnack (also
called Li-Yau-Hamilton) expression for shrinking solitons, $X$ is
$(1,0)$ component of  the tangent vector of the minimizing ${\Cal
L}$-geodesics. Notice that $\langle\cdot, \cdot\rangle$ is the
Hermitian product with respect to the K\"ahler metric.  We also
have that
$$
|\nabla \ell|^2+\ell_\tau=-\frac{1}{2\tau^{\frac{3}{2}}}K.\tag
1.12
$$
and
$$
\ell_\tau={\Cal
R}-\frac{\ell}{\tau}+\frac{1}{2\tau^{\frac{3}{2}}}K. \tag 1.13
$$
Applying H.-D. Cao's \cite{C2} (in stead of Hamilton's) trace
differential Harnack (also called Li-Yau-Hamilton inequality) for
 ancient solutions to K\"ahler-Ricci flow,
$$
-{\Cal R}_{\tau}-2\langle\nabla {\Cal R}, X\rangle-2\langle X,
\nabla {\Cal R}\rangle +4Ric(X, X)\ge 0 \tag 1.4
$$
we have  that
$$
H(X)\ge -\frac{1}{\eta}{\Cal R}.
$$
Therefore we have that
$$
K\ge -\int_0^\tau \sqrt{\eta}{\Cal R}\, d \eta
\ge-2\sqrt{\tau}\ell.
$$
Applying the above lower bound to (1.11)--(1.13) we have that
$$
|\nabla \ell|^2+{\Cal R}\le \frac{3}{\tau} \ell, \tag 1.15
$$
$$
|\nabla \ell|^2+\ell_\tau \le \frac{\ell}{\tau} \tag 1.16
$$
and
$$
\ell_\tau \ge {\Cal R}-\frac{2}{\tau} \ell\ge -\frac{2}{\tau} \ell
. \tag 1.17
$$
From (1.15), we have that
$$
|\nabla \ell^{\frac{1}{2}}|^2\le \frac{3}{4\tau},
$$
which implies that
$$
\ell^{\frac{1}{2}}(y, \tau)\le
\sqrt{\frac{n}{2}}+d_{\tau}(x(\tau), y)\sqrt{\frac{3}{4\tau}}.
\tag 1.18
$$
Now we can deduce the following results on the reduced distance
$\ell (y, \tau)$. \proclaim{Corollary 1.1} For $y\in
B_{\tau}(x(\tau), \sqrt{\frac{\tau}{\e}})$,
$$
\ell(y, \tau)\le
\lf(\sqrt{\frac{n}{2}}+\sqrt{\frac{3}{4\e}}\ri)^2. \tag 1.19
$$
Hence
$$
\ell(y, \eta)\le
4\lf(\sqrt{\frac{n}{2}}+\sqrt{\frac{3}{4\e}}\ri)^2. \tag 1.20
$$
for all $2\tau\ge \eta\ge \frac{\tau}{2}$ and $y\in
B_{\tau}(x(\tau), \sqrt{\frac{\tau}{\e}})$. And
$$
\tau {\Cal R}(y, \eta) \le
12\lf(\sqrt{\frac{n}{2}}+\sqrt{\frac{3}{4\e}}\ri)^2. \tag 1.21
$$
Moreover, on $ B_{\tau}(x(\tau), \sqrt{\frac{\tau}{\e}})\times
(\delta \tau, \frac{1}{\delta} \tau)$ one has that
$$
\ell(y, \eta)\le
\lf(\frac{1}{\delta}(\sqrt{\frac{n}{2}}+\sqrt{\frac{3}{4\e}})\ri)^2.
\tag 1.20'
$$
and
$$
\tau {\Cal R}(y, \eta) \le
3\lf(\frac{1}{\delta}(\sqrt{\frac{n}{2}}+\sqrt{\frac{3}{4\e}})\ri)^2.
\tag 1.21'
$$
\endproclaim
The following corollary gives relation between the
$\kappa$-constant in the definition of $\kappa$-noncollapsing and
the lower bound of the reduced volume. The converse is also true,
even without any curvature sign assumptions. Please see \cite{CLN}
and author's Dec, 2003 lecture at AIM for details.

\proclaim{Corollary 1.2} There exists a constant $C_1(n)>0$ such
that
$$
\tilde{V}(\tau)\ge e^{-C_1(n)}\kappa. \tag 1.22
$$
This in particular implies that $\lim_{\tau \to \infty}
\tilde{V}(\tau)>0.$
\endproclaim
\demo{Proof} Apply the (1.21)  to $B_\tau(x(\tau),\tau)$ to see
that it satisfies the curvature bound assumption on the ball in
the non-collapsing assumption. Then the result follows from the
estimate (1.20) and the definition of the reduced volume $\tilde
V(\tau)=\int_M \frac{e^{-\ell(y\, \tau)}}{\tau^{\frac{n}{2}}}\,
d\mu_\tau$.
\enddemo

Let $\tg^\tau(s)=\frac{1}{\tau}\bg(s\tau)$. Then
$\tilde{B}_{1}(x(\tau), \sqrt{\frac{1}{\e}})=B_{\tau}(x(\tau),
\sqrt{\frac{\tau}{\e}})$. The (1.21') and the
$\kappa$-non-collapsing assumption implies, by Hamilton's
compactness result \cite{H3}, that
$$\lf(\tilde{B}_{1}(x(\tau),
\sqrt{\frac{1}{\e}})\times(\delta, \frac{1}{\delta}),
\tg^{\tau}(s)\ri)\to \lf(B^{\infty}(x_{\infty},
\sqrt{\frac{1}{\e}})\times(\delta, \frac{1}{\delta}),
g^{\infty}(s)\ri)$$ as solutions to K\"ahler-Ricci flow. This can
 be extended to an  ancient solution  $(M_{\infty},
g_{\infty}(\tau))$. The estimates (1.15)--(1.17) ensure that
$\ell(y, s)$ the reduced distance with respect to re-scaled metric
survives under the limit and converges to a function
$\ell_{\infty}(y, s)$. Let
$$
V_{\infty}(s)=\int_{M_{\infty}}\frac{e^{-\ell_{\infty}(y,s)}}{s^{\frac{n}{2}}}\,
d\mu_{s}.
$$
It was claimed by Perelman in \cite{P1, 11.2} that
 \proclaim{Claim 2}
 $$
\lim_{\tau \to \infty}\tilde{V}(\tau)=V_{\infty}(s). \tag 1.23
$$ In particular, one has  that $V_{\infty}(s)$ is a constant and $(M_\infty, g^{\infty})$
is a non-flat gradient shrinking soliton.
\endproclaim
The claim follows if one can show that
$$
\int_{(B_{\tau}(x(\tau),
\sqrt{\frac{\tau}{\e}}))^c}\frac{e^{-\ell(y,
\tau)}}{\tau^{\frac{n}{2}}}\, d\mu_{\tau} \le C(\e) \tag 1.24
$$
with $C(\e)\to 0$ as $\e\to 0$. This  can be proved easily if we
can get an `effective' lower bound estimate of $\ell(y, \tau)$ in
terms of $\frac{d^2_\tau(x_0, y)}{\tau}$. This last point   was
proved in \cite{Ye, Lemma 2.2}. In fact it was proved that there
exists positive constant $C=C(n)$ such that for any $y, z\in M$.
$$
-\ell(z, \tau)-1+C\frac{d_{\tau}^2(z, y)}{\tau}\le \ell(y, \tau).
$$
One should consult notes \cite{KL, Ye} for more detailed
exposition on the proof of Claim 2.

\demo{Proof of Theorem 2} Assume that $(M, g(t))$ is an ancient
solution (non-flat) defined on $M\times (\infty, 0]$. If $M$ is
compact, there is nothing to prove. So we  assume that $M$ is
non-compact. We prove the theorem by the contradiction. So we
assume that ${\Cal V}(g(t_0))>0$ for some $t_0$. By passing to it
universal cover we can also assume that $M$ is simply-connected.
 Then by part (ii) of
Theorem 4.2 of \cite{NT2} again (see also Corollary 4.1 of cited
paper), we have that the scalar curvature has the average decay
(1.10). Now apply Theorem 2.2 of \cite{NT1} (and its proof for the
time before $t_0$) we conclude that for all $t$, ${\Cal
V}(g(t))={\Cal V}(g(t_0))>0$. This in particular implies that $(M,
g(t))$ is $\kappa$-non-collapsed in all scales (since the volume
is non-collapsed even without assuming the curvature bound).

Now we apply the above blow-down procedure of Perelman to obtain
the limit $(M_\infty, g^{\infty})$, which is a gradient shrinking
soliton. Since $(M, g(t))$ is assumed to have  nonnegative
bisectional curvature, the limit $(M_\infty, g^{\infty})$ also has
 nonnegative bisectional curvature. By the  definition of
$(M_{\infty}, g^{\infty})$, it is clear that the corresponding
asymptotic volume ratio ${\Cal V }(M_\infty, g^{\infty})\ge {\Cal
V}(g(t_0))>0$. On the other hand, since $(M_\infty, g^{\infty})$
is non-flat, its scalar curvature must be positive by the strong
maximum principle. Now  we can apply Theorem 3 (part (ii)) to the
universal cover of $M_{\infty}$ and conclude that $M_\infty$ can
not have maximum volume growth. This contradicts the fact that
${\Cal V }(M_\infty, g^{\infty})\ge {\Cal V}(g(t_0))>0$.  The
contradiction proves the theorem.

\enddemo

\subheading{\S2 Applications to K\"ahler manifolds and
K\"ahler-Ricci flow}\vskip .2cm

Theorem 1 has nice applications to Ricci flow as shown in
\cite{P1}. Following \cite{P1}, we can derive  the  compactness
result on ancient solutions to K\"ahler-Ricci flow out of Theorem
2.

We call an ancient non-flat solution $(M, g(t))$ defined on
$M\times (-\infty, 0]$ a {\it $\kappa$-solution} to K\"ahler-Ricci
flow, if $(M, g(t))$ has nonnegative bisectional curvature,
satisfying the trace differential Harnack inequality
$$
{\Cal R}_t+\langle \nabla {\Cal R}, X\rangle +\langle X, \nabla
{\Cal R}\rangle + Ric(X, X)\ge 0 \tag 2.1
$$
for any $(1,0)$ vector field $X$,  and  $(M, g(t))$ is
$\kappa$-non-collapsed on all scales for some fixed $\kappa>0$.
Recall that $\kappa$-non-collapsed on all scales means that for
any time $t$ and $x_0\in M$ if for all $y\in B_t(x_0, r)$, ${\Cal
R}(y, t)\le \frac{1}{r^2}$, then $Vol_t(B_t(x_0, r))\ge \kappa
r^{2m}$. Notice that we do not require ${\Cal R}$ being bounded.
We formulate in such way to be able to make the result hold for
any dimension.

\proclaim{Theorem 2.1} The set of $\kappa$-solutions to
K\"ahler-Ricci flow  is compact modulo scaling.
\endproclaim

The proof of the result follows the same line of the argument as
Theorem 11.7 of \cite{P1}. The argument in \cite{P1} is quite
robust. The key components of the argument in \cite{P1} are Shi's
local derivative estimate, trace differential Harnack and the
following consequence of Theorem 2.

\proclaim{Corollary 2.1} For every $\omega>0$ there exist
$B=B(\omega)<\infty ,
  C=C(\omega)<\infty , \tau_0=\tau_0(\omega)>0,$ with the following properties.
  \par (a) Suppose we have a (not necessarily complete) solution
  $g(t)$ to the K\"ahler-Ricci flow, defined on $M\times [t_0,0],$ so
  that at time $t=0$ the metric ball $B_0(x_0,r_0)$ is compactly
  contained in $M.$ Suppose that at each time $t, t_0\le t\le 0,$
  the metric $g(t)$ has nonnegative bisectional curvature, and
  $Vol_t (B_t(x_0,r_0))\ge \omega r_0^n.$ Then we have an estimate $R(x,t)\le
  Cr_0^{-2}+B(t-t_0)^{-1}$ whenever $d_t(x,x_0)\le
  \frac{1}{4}r_0.$
  \par (b) If, rather than assuming a lower bound
  on volume for all $t,$ we
 assume it only for $t=0,$ then the same conclusion holds with
 $-\tau_0r_0^2$ in place of $t_0,$ provided that $-t_0\ge
 \tau_0r_0^2.$
\endproclaim

The corollary above is exactly the same as Corollary 11.6 of
\cite{P1} with only the Ricci flow being replaced by
K\"ahler-Ricci flow and curvature operator being replaced by
bisectional curvature. The proof is the also the same if one
replaces   Theorem 1 of Perelman, whereever it is needed,   by
Theorem 2 of  this paper. The robust scaling argument in the proof
of Corollary 11.6 in \cite{P1} also resembles, to some degree in
the spirit,  various curvature estimates proved by certain scaling
argument in the study of the mean curvature flow and other PDEs.
See for example, \cite{S1, E, Wh}, as well as Simon's proof on
Schauder's estimates in \cite{S2}. The detailed exposition on the
proof of the above Theorem 2,1 and Corollary 2.1 can be found in
\cite{CLN, KL, STW} as well as author's AIM lecture notes.

In view of the compactness result above and the general tensor
maximum principle proved in \cite{NT2}, we conjecture that H.-D.
Cao's matrix Li-Yau-Hamilton estimate holds on any complete
K\"ahler manifolds with non-negative bisectional curvature. No
assumption on curvature {\it being bounded} is needed. If
confirmed, one has Theorem 2.1 for all $\kappa$ non-collapsed
ancient solutions with non-negative bisectional curvature. As in
\cite{P1}, the following gradient estimate is a consequence of
Theorem 2.1.

\proclaim{Corollary 2.2} There exists $C=C(\kappa, m)$ such that
for the $\kappa$-solution $(M, g(t))$ we have that
$$
|\nabla {\Cal R}|(x,t)\le C {\Cal R}^{\frac{3}{2}}(x,t), \quad
\quad |{\Cal R}_t|(x,t)\le C {\Cal R}^2(x,t).
$$
\endproclaim

Note that Theorem 2  also holds  for ancient solutions with {\it
nonnegative bisectional curvature and differential Harnack (2.1).}
 Also notice that the {\it $\kappa$-solution} here has different meaning
from \cite{P1. Sectiona 11}. The {\it bounded $\kappa$-soltion}
defined last section corresponds to the $\kappa$-solution in
\cite{P1}.

\demo{Proof of Corollary 1} By Shi's short time existence result
we know that the solution exists until the curvature blows up. But
by Corollary 2.1 above (see also \cite{P1, 11.5 and 11.6}) one has
the estimate
$$
{\Cal R}(x,t)\le \frac{C_1}{t+1}, \tag 2.2
$$
 for some $C_1=C_1(n,{\Cal
V}(g(0)))>0$ Therefore, the solution exists for all time and is of
Type III. Namely there is no slowly forming singularity at
infinity. The result can also be shown using Hamilton's blow-up
argument in \cite{H5}. In order to get the curvature decay
estimate (0.2) we first apply Theorem 2.1 of \cite{NT1} to
conclude that there exists $C_2=C_2(m)$ with
$$
\int_0^{\sqrt{t}}sk(x, s)\, ds \le -C_2F(x, t)
$$
where $k(x,r)=\aint_{B_x(r)}{\Cal R}(y, 0)\, d\mu$, with respect
to the initial metric and
$$
F(x,t)=\log\left[\frac{\det(g_{\abb}(x,t))}{\det(g_{\abb}(x,0))}\right].$$
(The above  estimate follows easily from that facts $-F(x,t)\ge 0$
and  $\left(\D_0-\frac{\p}{\p t}\right)(-F)(x,t)\le -{\Cal
R}(x,0)$, where $\D_0$ is the Laplace operator with respect to
$g(x,0)$. This is indeed the proof on page 125 of \cite{NT1}.)
Since $-\frac{\p}{\p t}F={\Cal R}(x, t)\le \frac{C_1}{1+t}$, for
$1\ll t$, one has that
$$
\int_0^{r}sk(x, s)\, ds\le C_3\log (r+2),\tag 2.3
$$
for some $C_3=C_3(m, {\Cal V}(g(0)))$. (This implies what proved
in \cite{CZ4} under the stronger assumption on the nonnegativity
of curvature operator.)

Using the curvature decay estimate (2.2) together with the fact
that the asymptotic volume ratio ${\Cal V}(t)$ is a constant
function of $t$, one can apply the local injectivity radius
estimate of Cheeger-Gromov-Taylor \cite{CGT, Theorem 4.3} (see
also \cite{CLY} for earlier similar works)  to conclude that the
injectivity radius of $(M, g(t))$ has the size of $C\sqrt{t}$,
where $C$ is a constant independent of $t$. This
 implies that $M$ can be exhausted by open subdomains which
diffeomorphic to Euclidean balls. From this one can conclude the
topological type of $M$ from by-now standard result from topology.
The above is the observation in \cite{CZ3}, which follows the
earlier construction in \cite{Sh4, Section 9} and \cite{H6}. Note
that the Steinness of $M$ has been proved for any complete
K\"ahler manifolds with the maximum volume growth and nonnegative
bisectional curvature in \cite{NT2, Corollary 4.2} (Please see
also \cite{WZ} for the even easier positive case.) One then can
adapt the construction of \cite{Sh4, Section 9} to conclude that
$M$ is biholomorphic to a pseudoconvex domain in $\C^m$.

To obtain (0.2), by \cite{NT2, Theorem 6.1}, we first  solve the
Poincar\'e-Lelong equation to obtain  $u$ such that $\p_{\a}
\bar{\p}_{\bar{\beta}}u=R_{\abb}$, and $u$ is at most of
logarithmic growth. By considering its heat equation deformation
and adding  a function of the form $\log (|z|^2+1)$ in the case
$M$  splits some factors of $\C$, noticing that $M$ is
diffeomorphic to $\R^n$, we can obtain a strictly plurisubharmonic
function of logarithmic growth on M.
 More precisely, let $v(x,t)$ be the heat equation deformation of $u(x)$
(without K\"ahler-Ricci flow). By \cite{NT2}, $v(x,t)$ has the
 same growth as $u(x)$ and $M$ splits as $M=M_1\times M_2$, where on $M_1$ the complex
Hessian $v_{\abb}(x,t)$ is positive definite and on $M_2$,
$v_{\abb}(x,t)\equiv 0$. By the result of Cheng-Yau we know that
$v(x,t)$ must be a constant on $M_2$, which then implies that
$u(x)$ is a constant on $M_2$. This implies that $M_2=\C^k$ for
some $0\le k\le m$. On $M_2$ one can construct a strictly
plurisubharmonic function of logarithmic growth by taking
$\phi(z)=\log (|z|^2+1)$. Adding this to $v(x,t)$ we have a
strictly plurisubharmonic function of logarithmic growth on  $M$.
Now (0.2) follows from the proof of Corollary 3.2 in  \cite{N2}.
\enddemo

\demo{Proof of Corollary 2} By the above proof of Corollary 1  we
know that $M$ supports a strictly plurisubharmonic function of
logarithmic growth. The conclusion on the transcendence degree of
rational function field  follows from the $L^2$-estimate of
$\bar{\p}$ and the dimension estimate proved in Theorem 3.1 of
\cite{N2}. See, for example, \cite{N2, Theorem 5.2} for the
constructions of holomorphic functions of polynomial growth,
forming a local holomorphic coordinate for any given point in $M$.
The construction via the well-known H\"ormander's $L^2$-estimates
provides the existence of plenty holomorphic functions of
polynomial growth, which provides a lower bound on the
transcendence degree. The multiplicity estimates in Theorem 3.1 of
\cite{N2} gives the upper bound on the transcendental elements of
the rational functions. The assertion that $M$ is an affine
quasi-algebraic variety follows from the construction in
\cite{M1}. See also \cite{De} as well as \cite{CTZ}. Notice that
we now also have uniform multiplicity estimates, thanks to the new
monotonicity formula and Theorem 3.1 in \cite{N2}, which
simplifies the steps of \cite{M1} quite a bit. Here by an affine
quasi-algebraic variety we mean an affine algebraic variety with
some codimension one algebraic varieties removed.

The part $M$ is biholomorphic to $\C^2$ also follows as in
\cite{M1, page 256} by appealing to the result of Ramanujam that
any quasi-projective surface homeomorphic to $\R^4$ is
biholomorphic to $\C^2$.
\enddemo

\proclaim{Remark 2.1} The fact that maximum volume growth implies
the ampleness of holomorphic functions of polynomial growth was
conjectured in \cite{N2}. It has been shown in \cite{NT2} that the
positivity  of Ricci curvature together with some average
curvature decay assumption also implies the same result on
holomorphic function of polynomial growth. In fact, this also
implies that $M$ is an affine quasi-algebraic variety as in
Corollary 2.

The existence of harmonic functions of polynomial growth was
obtained in \cite{D} under the assumption of nonnegative Ricci
curvature, maximum volume growth and the uniqueness of the tangent
cone at infinity. The result here seems to indicate that the
assumption on the  uniqueness of the tangent cone may not be
needed.

After the completion of current paper, H.-D. Cao informed the
author that the case of $m=2$ of Corollary 2,  assuming the
maximum volume growth and  {\it nonnegativity of bisectional
curvature}, was also obtained by B.-L. Chen earlier.
\endproclaim

\subheading{\S3 Applications to Ricci flow}\vskip .2cm

 The proof of Proposition 1.1   can be used in some
other situations. For example we can prove the following results.

\proclaim{Proposition 3.1} \roster \item"{(i)}" Let $(M, g)$ be a
compete steady gradient soliton. Assume that the Ricci curvature
is pinched in the sense that $R_{ij}\ge \e {\Cal R}g_{ij}$, for
some $\e>0$ with the scaler curvature ${\Cal R}(x)>0$. Then there
exist $C>0$, $a>0$ such that
$$
{\Cal R}(x)\le C\exp(-a(r(x)+1)). \tag 3.1
$$
Here $r(x)$ is the distance function to some fixed point in $M$.

\item"{(ii)}" Let $(M, g)$ be a compete expanding gradient
soliton. Assume that the Ricci curvature is pinched as above. Then
{\rm (3.1)} holds.
\endroster
\endproclaim
\proclaim{Corollary 3.1} Assume that $n\ge 3$. \roster
\item"{(i)}" There is no steady gradient soliton with pinched
Ricci curvature as in Proposition 3.1;

\item"{(ii)}" There is no expanding  gradient soliton with pinched
Ricci curvature as in Proposition 3.1 and nonnegative sectional
curvature.
\endroster
In particular, any complete three manifolds with pinched Ricci
curvature must be compact, therefore spherical. Same result holds
for any higher dimensional complete Riemannian manifolds with
pinched curvature operator in the sense that $$ |\oR m|^2=
|W|^2+|V|^2\le \delta_n (1-\epsilon)^2|U|^2=\delta_n
(1-\epsilon)^2\frac{2}{n(n-1)}{\Cal R}^2 \tag 3.2
$$
where $\e>0$, $\delta_3>0, \delta_4=\frac{1}{5},
\delta_5=\frac{1}{10}$ and $\delta_n =\frac{2}{(n-2)(n+1)}$, and
$W$, $V$ and $U$ denote the Weyl curvature tensor, traceless Ricci
part and the scalar curvature part, respectively, according to the
curvature operator decomposition in \cite{Hu}.
\endproclaim
\demo{Proof} The first part follows from Theorem 20.2 of \cite{H5}
and (i) of Proposition 3.1. Notice that the proof there works
under the weaker assumption $R_{ij}>0$.

For part (ii), one just need to recall the gap theorem of
Greene-Wu \cite{GW} (see also \cite{PT}) asserting that a
simply-connected complete Riemannian manifold with nonnegative
sectional curvature and (3.1) must be flat, noticing that under
the assumption of expanding gradient soliton and nonnegativity of
the Ricci curvature, $M$ is diffeomorphic to $\R^n$ (which in
particular implies that the new
 expanders examples constructed in \cite{FIK} can not have
 nonnegative Ricci curvature), since $f$ is a strictly convex function with only one
  critical point.

  By (i) and (ii) and the discussion above, we conclude that curvature pinched manifolds
  must be compact. This is the main result proved in \cite{CZ1}.
  Now  the last  claim in the corollary just restates the fundamental results
  of Hamilton \cite{H1} and Huisken \cite{Hu}.
\enddemo

Note that when $n=2$, (3.1) is automatic and the examples of
Hamilton's cigar and the expanding soliton exhibited in \cite{CLN}
show that the exponential decay proved in Proposition 3.1 is
sharp.

\proclaim{Remark 3.1} We speculate that there is no complete
noncompact Riemannian manifold with the pinched Ricci as in
Proposition 3.1 (we found out later that, according to Ben Chow,
this was asked in case $n=3$ by Hamilton earlier). If true, this
is really a new Bonnet-Meyers type theorem since (3.2) is too
strong to allow other topology. But we do not have any workable
scheme to prove such general result at this moment. The details on
the proof of results in this section will appear somewhere else.

\endproclaim

\enddocument